\documentstyle[11pt]{article}
 \newlength{\standardunitlength}
\newtheorem{cor}{Corollary} \newtheorem{lemma}{Lemma}
\newtheorem{theorem}{Theorem} 
\newenvironment{proof}{\noindent {\sc Proof:}}{$\Box$ \vspace{2 ex}}
\begin{document}

\begin{center} An Application of the Whitehouse Module to Riffle
Shuffles Followed by a Cut \end{center}

\begin{center}
By Jason Fulman
\end{center}

\begin{center}
Stanford University
\end{center}

\begin{center}
Department of Mathematics
\end{center}

\begin{center}
Building 380, MC 2125
\end{center}

\begin{center}
Stanford, CA 94305
\end{center}

\begin{center}
email:fulman@math.stanford.edu
\end{center}

\begin{center}
http://math.stanford.edu/$\sim$fulman
\end{center}

\newpage \begin{abstract} Using representation theoretic work on the
Whitehouse module, a formula is obtained for the cycle structure of a
riffle shuffle followed by a cut. This result will be merged with the
paper \cite{F6}. \end{abstract}

\section{Introduction}

	In an effort to study the way real people shuffle cards, Bayer
and Diaconis \cite{BD} performed a definitive analysis of the
Gilbert-Shannon-Reeds model of riffle shuffling. For an integer $k
\geq 1$, a $k$-shuffle can be described as follows. Given a deck of
$n$ cards, one cuts it into $k$ piles with probability of pile sizes
$j_1,\cdots,j_k$ given by $\frac{{n \choose
j_1,\cdots,j_k}}{k^n}$. Then cards are dropped from the packets with
probability proportional to the pile size at a given time (thus if the
current pile sizes are $A_1,\cdots,A_k$, the next card is dropped from
pile $i$ with probability $\frac{A_i}{A_1+\cdots+A_k}$). It was proved
in \cite{BD} that $\frac{3}{2}log_2n$ shuffles are necessary and
suffice for a $2$-shuffle to achieve randomness (the paper \cite{A}
had established this result asymptotically in $n$). It was proved in
\cite{DMP} that if $k=q$ is a prime power, then the chance that a
permutation distributed as a $q$-shuffle has $n_i$ $i$-cycles is equal
to the probability that a uniformly chosen monic degree $n$ polynomial
over the field $F_q$ factors into $n_i$ irreducible polynomials of
degree $i$.

	A very natural question is to study the effects of cuts on the
results of the previous paragraph. For example, it is shown in
\cite{F6} that performing the process of ``a riffle shuffle followed
by a cut at a uniform position'' also gets random in
$\frac{3}{2}log_2n$ steps. This can be contrasted with a result of
Diaconis \cite{D}, who proves that although shuffling by doing random
tranpositions gets random in $\frac{1}{2}n log(n)$ steps, the use of
cuts at each stage drops the convergence time to $\frac{3}{8} n
log(n)$ steps. The main result of this note is that if $k=q$ is a
prime power, then the chance that a permutation distributed as a
$q$-shuffle has $n_i$ $i$-cycles is equal to the probability that a
uniformly chosen monic degree $n$ polynomial over the field $F_q$ {\it
with non-zero constant term} factors into $n_i$ irreducible
polynomials of degree $i$. The result is proved by showing it to be
completely equivalent to representation theoretic results of
Whitehouse \cite{W}.

	Before jumping into the proof, we remark that the theory of
riffle shuffling appears in numerous parts of mathematics. Among these
are:

\begin{enumerate}

\item Cyclic and Hochschild homology \cite{H},\cite{GS},\cite{L}
\item Hopf algebras, Poincar\'e-Birkhoff-Witt theorem (Chapter 3.8 of \cite{SS})
\item Representation theory of the symmetric group \cite{Sta}
\item Dynamical systems \cite{BD},\cite{La1},\cite{La2},\cite{F5}
\item Free Lie algebras \cite{Ga}
\item Random matrices \cite{Sta}
\item Algebraic number theory (speculative) \cite{F3},\cite{F5}
\item Potential theory \cite{F7}

\end{enumerate} A survey paper describing these connections, to be
titled ``Riffle shuffling: a unifying theme'' is in preparation.

	In the past few years interesting combinatorial
generalizations of riffle shuffling have emerged. Roughly, they can be
classified as

\begin{enumerate}
\item Biased riffle shuffles \cite{DFP},\cite{F1},\cite{Sta}
\item Other Coxeter groups \cite{BD}, \cite{BB}, \cite{F2}
\item Hyperplane arrangements \cite{BHR}, \cite{F2}
\item Affine shuffles \cite{Ce1},\cite{Ce2},\cite{F5},\cite{F6}
\item Riffle shuffles with cuts \cite{BD}, \cite{Ce3}, \cite{F6}
\end{enumerate}

	The point is that since riffle shuffles are related to so many
parts of mathematics, these generalizations should be interesting
too. In particular, as is clear from \cite{DMP} and the papers of the
author just cited (see also \cite{Sta} for connections with
quasi-symmetric functions and extensions to infinite support) for
these generalizations the induced distribution on conjugacy classes
seems to have lovely properties. This note gives further support to
that philosophy.

	We remark that the Whitehouse module also appears in
interesting mathematical contexts (\cite{HS},\cite{W2},
\cite{LS}). Richard Stanley's MIT website contains transparencies from
an illuminating talk about the Whitehouse module.

	The structure of this note is as follows: Section \ref{main}
gives the main result, and Section \ref{open} suggests two open
problems.

\section{Main Result} \label{main}

	To begin some notation is necessary. Recall that an element
$w$ of $S_n$ is said to have a descent at position $i$ (with $1 \leq i
\leq n-1$) if $w(i)>w(i+1)$, and a cyclic descent at position $n$ if
$w(n)>w(1)$. One lets $d(w)$ be the number of descents of $w$ and
defines $cd(w)$ to be $d(w)$ if $w$ has no cyclic descent at $n$, and
to be $d(w)+1$ if $w$ has a cyclic descent at $n$. Thus $cd(w)$ can be
thought of as the total number of descents of $w$, viewed cyclically.

	Now we use representation theory to obtain a formula for the
cycle structure of a riffle shuffle followed by a cut. To begin we
recall the following result which gives a formula for the chance of a
permutation $w$ after a $k$-riffle shuffle followed by a cut.

\begin{theorem} (\cite{F6}) \label{correct} The chance of obtaining a
permutation $w$ after a $k$-riffle shuffle followed by a cut is \[
\frac{1}{nk^{n-1}} {n+k-cd(w^{-1})-1 \choose n-1}.\] \end{theorem}

	It is useful to recall the notion of a cycle index associated
to a character of the symmetric group. Letting $n_i(w)$ be the number
of $i$-cycles of a permutation $w$ and $N$ be a subgroup of $S_n$, one
defines $Z_N(\chi)$ as

\[ Z_N(\chi) = \frac{1}{|N|} \sum_{w \in N} \chi(w) \prod_i
a_i^{n_i(w)}.\] The cycle index stores complete information about the
character $\chi$. For a proof of the following attractive property of
cycle indices, see \cite{Fe}.

\begin{lemma} \label{Walter} Let $N$ be a subgroup of $S_n$ and $\chi$
a class function on $N$. Then 

\[ Z_{S_n}(Ind_N^{S_n}(\chi)) = Z_N(\chi).\]
\end{lemma}

	Next, recall that an idempotent $e$ of the group algebra of a
finite group $G$ defines a character $\chi$ for the action of $G$ on
the left ideal $KGe$ of the group algebra of $G$ over a field $K$ of
characteristic zero. For a proof of Lemma \ref{bridge}, which will
serve as a bridge between representation theory and computing measures
over conjugacy classes, see \cite{H}. For its statement, let $e<w>$ be
the coefficient of $w$ in the idempotent $e$.

\begin{lemma} \label{bridge} Let $C$ be a conjugacy class of the
finite group $G$, and let $\chi$ be the character associated to the idempotent
$e$. Then

\[ \frac{1}{|G|} \sum_{w \in C} \chi(w) = \sum_{w \in C} e<w>.\]
\end{lemma}

	 It is also convenient to define

\[ Z_{S_n}(e) = \sum_{w \in S_n} e<w> \prod_i a_i^{n_i(w)},\] which
makes sense for any element $e$ of the group algebra. Note that one
does not divide by the order of the group. When $e$ is idempotent and
$\chi$ is the associated character, Lemma \ref{bridge} can be
rephrased as

\[ Z_{S_n}(\chi) = Z_{S_n}(e).\]

	To proceed recall the Eulerian idempotents $e_n^j$,
$j=1,\cdots,n$ in the group algebra $QS_n$ of the symmetric group over
the rationals. These can be defined \cite{GS} as follows. Let
$s_{i,n-i}=\sum w$ where the sum is over all ${n \choose i}$
permutations $w$ such that $w(1)<\cdots<w(i)$, $w(i+1)<\cdots<w(n)$
and let $s_n=\sum_{i=1}^{n-1} s_{i,n-i}$. Letting $\mu_j=2^j-2$, the
$e_n^j$ are defined as

\[ e_n^j = \prod_{i \neq j} \frac{s_n-\mu_i}{(\mu_j-\mu_i)}.\] They
are orthogonal idempotents which sum to the identity.

	The following result, which we shall need, is due to
Hanlon. The symbol $\mu$ denotes the Moebius function of elementary
number theory.

\begin{theorem} (\cite{H}) \label{Hanl}

\[ 1+\sum_{n=1}^{\infty} \sum_{i=1}^n k^i Z_{S_n}(e_n^i) = \prod_{i
\geq 1} (1-a_i)^{-(1/i) \sum_{d|i} \mu(d) k^{i/d}}.\] \end{theorem}

\begin{theorem} (\cite{Ga}) \label{Gars}

\[ \sum_{i=1}^n k^i e_n^i = \sum_{w \in S_n} {n+k-d(w)-1 \choose n}
w.\] \end{theorem}

{\bf Remark:} Combining Lemma \ref{bridge} and Theorem \ref{Gars}, one
sees that the formula for the cycle structure of a riffle shuffle
\cite{DMP} and Theorem \ref{Hanl} imply each other. It is interesting
that both proofs used a bijection of Gessel and Reutenauer \cite{G}.

	To continue, we let $\overline{e_n^j}$ denote the idempotent
obtained by multiplying the coefficient of $w$ in $e_n^j$ by
$sgn(w)$. Let $\lambda_{n+1}$ be the $n+1$ cycle $(1 \ 2 \cdots n+1)$
and $\Lambda_{n+1}=\frac{1}{n+1} \sum_{i=0}^n (sgn \lambda_{n+1}^i)
\lambda_{n+1}^i$. Viewing $\overline{e_n^j}$ as in the group algebra
of $S_{n+1}$, Whitehouse \cite{W} proves that for $j=1,\cdots,n$ the
element $\Lambda_{n+1} \overline{e_n^j}$ is an idempotent in the group
algebra $QS_{n+1}$, which we denote by $f_{n+1}^j$. Whitehouse's main
result is the following:

\begin{theorem} (\cite{W}) \label{white} Let $F_{n+1}^j,\overline{E_{n}^j}$
be the irreducible modules corresponding to the idempotents
$f_{n+1}^j$ and $\overline{e_n^j}$. Then

\[ F_{n+1}^j \oplus \bigoplus_{i=1}^j \overline{E_{n+1}^i} =
\bigoplus_{i=1}^j Ind_{S_n}^{S_{n+1}} \overline{E_{n}^i}.\] \end{theorem}

	As final preparation for the main result of this section, we
link the idempotent $\Lambda_{n+1} \overline{e_n^j}$ with riffle shuffles
followed by a cut.

\begin{lemma} \label{form} The coefficient of $w$ in $\sum_{j=1}^n k^j
\Lambda_{n+1} \overline{e_n^j}$ is $sgn(w) \frac{1}{n+1} {k+n-cd(w)
\choose n}$. \end{lemma}

\begin{proof} Given Theorem \ref{Gars}, this is an elementary
combinatorial verification. \end{proof}

	Theorem \ref{cyc} now derives the cycle structure of a
permutation distributed as a shuffle followed by a cut. So as to
simplify the generating functions, recall that $\sum_{d|i} \mu(d)$
vanishes unless $i=1$.

\begin{theorem} \label{cyc}
\begin{eqnarray*}
&& 1+\sum_{n \geq 1} \sum_{w \in S_{n+1}}
\frac{1}{(n+1)k^{n+1}} {n+k-cd(w) \choose n} \prod_i a_i^{n_i(w)}\\
& = & 1-\frac{1}{k-1} -\frac{a_1}{k} + \frac{1}{k-1} \prod_{i \geq 1}
(1-\frac{a_i}{k^i})^{-1/i \sum_{d|i} \mu(d) (k^{i/d}-1)}.
\end{eqnarray*} If $k=q$ is the size of a finite field, this says that the cycle type of a
permutation distributed as a shuffle followed by the cut has the same
law as the factorization type of a monic degree $n$ polynomial over
$F_q$ with non-vanishing constant term.  \end{theorem}

\begin{proof} Replacing $a_i$ by $a_i k^i (-1)^{i+1}$, it is enough to
show that

\begin{eqnarray*}
& & 1+\sum_{n \geq 1} \sum_{w \in S_{n+1}} sgn(w) \frac{1}{(n+1)}
{n+k-cd(w) \choose n} \prod_i a_i^{n_i(w)}\\
& = &1-\frac{1}{k-1} -a_1 +
\frac{1}{k-1} \prod_{i \geq 1} (1-(-1)^{i+1}a_i)^{-1/i \sum_{d|i}
\mu(d) (k^{i/d}-1)}.
\end{eqnarray*}

	Using Lemmas \ref{Walter}, \ref{bridge}, \ref{form} and
Theorem \ref{white}, one sees that

\begin{eqnarray*}
& &  1+\sum_{n \geq 1} \sum_{w \in S_{n+1}} sgn(w) \frac{1}{(n+1)}
{n+k-cd(w) \choose n} \prod_i a_i^{n_i(w)}\\
& =& 1+\sum_{n=1}^{\infty} \sum_{j=1}^n k^j Z_{S_{n+1}}(f_n^j)\\
& = & 1 + \sum_{n=1}^{\infty} \sum_{j=1}^n k^j \sum_{i=1}^j Z_{S_{n+1}}(Ind_{S_n}^{S_{n+1}}(\overline{e_n^i})) - \sum_{n=1}^{\infty} \sum_{j=1}^n k^j \sum_{i=1}^j Z_{S_{n+1}}(\overline{e_{n+1}^i})\\
& = & 1+a_1 \sum_{n=1}^{\infty} \sum_{i=1}^n Z_{S_n}(\overline{e_n^i})(\frac{k^{n+1}-k^i}{k-1}) - \sum_{n=1}^{\infty} \sum_{i=1}^n Z_{S_{n+1}}(\overline{e_{n+1}^i}) (\frac{k^{n+1}-k^i}{k-1})\\
& = &  1+a_1 k Z_{S_1}(\overline{e_1}) + \frac{a_1k-1}{k-1} \sum_{n=2}^{\infty} k^n \sum_{i=1}^n Z_{S_n}(\overline{e_n^i}) + \frac{1-a_1}{k-1} \sum_{n=2}^{\infty} \sum_{i=1}^n k^i Z_{S_n}(\overline{e_n^i}).
\end{eqnarray*} To simplify things further, recall that $\sum_{i=1}^n
Z_{S_n}(\overline{e_n^i})$ is $a_1^n$ since the $\overline{e_n^i}$'s sum to the
identity. The above then becomes

\[ 1-\frac{1}{k-1}-a_1+\frac{1-a_1}{k-1} (1+\sum_{n=1}^{\infty}
\sum_{i=1}^n k^i Z_{S_n}(\overline{e_n^i})),\] so the sought result follows
from Theorem \ref{Hanl}.
\end{proof}

	Before continuing, we observe that a combinatorial proof of
Theorem \ref{cyc} (which must exist) would give a new proof of Theorem
\ref{white}, by reversing the steps.

	Upon hearing about Theorem \ref{cyc}, Persi Diaconis
immediately asked for the expected number of fixed points after a
$k$-riffle shuffle followed by a cut, suggesting that it should be
smaller than for a $k$ riffle shuffle. Using the methods of Section 5
of \cite{DMP}, one can readily derive analogs of all of the results
there. As an illustrative example, Corollary \ref{fixpoint} shows that
the expected number of fixed points after a $k$-riffle shuffle
followed by a cut is the same as for a uniform permutation, namely 1
(the answer for $k$-riffle shuffles is $1+1/k+\cdots+1/k^{n-1}$). Two
other examples are worth mentioning and will be treated in Corollary
\ref{others}.

\begin{cor} \label{fixpoint} The expected number of fixed points after
a $k$-riffle shuffle followed by a cut is 1. \end{cor}

\begin{proof} The case $n=1$ is obvious. Multiplying $a_i$ by $u$ in
the statement of Theorem \ref{cyc} shows that

\begin{eqnarray*}
&& 1+\sum_{n \geq 1} \sum_{w \in S_{n+1}} u^{n+1}
\frac{1}{(n+1)k^{n+1}} {n+k-cd(w) \choose n} \prod_i a_i^{n_i(w)}\\
& = & 1-\frac{1}{k-1} -\frac{ua_1}{k} + \frac{1}{k-1} \prod_{i \geq 1}
(1-\frac{u^i a_i}{k^i})^{-1/i \sum_{d|i} \mu(d) (k^{i/d}-1)}.
\end{eqnarray*} To get the generating function in $u$ (for $n \neq 1$) for the expected number of fixed points in a riffle shuffle followed by a cut, one multiplies the right hand side by $k$, sets $a_2=a_3=\cdots=1$, differentiates with respect to $a_1$, and then sets $a_1=1$. Doing this yields the generating function \[ -u + u \prod_{i \geq 1} (1-\frac{u^i}{k^i})^{-1/i \sum_{d|i} \mu(d) k^{i/d}}.\] The result now follows from the identity

\[ \prod_{i \geq 1} (1-\frac{u^i}{k^i})^{-1/i \sum_{d|i} \mu(d)
k^{i/d}} = \frac{1}{1-u},\] which is equivalent to the assertion that
a monic degree $n$ polynomial over $F_q$ has a unique factorization
into irreducibles, since $1/i \sum_{d|i} \mu(d) k^{i/d}$ is the number
of irreducible polynomials of degree $i$ over the field $F_k$. \end{proof}

\begin{cor} \label{others} Fix $u$ with $0<u<1$. Let $N$ be chosen
from $\{0,1,2,\cdots\}$ according to the rule that $N=0$ with
probability $\frac{1-u}{1-u/k}$ and $N=n \geq 1$ with probability
$\frac{(k-1)(1-u)u^n}{k-u}$. Given $N$, let $w$ be the result of a
random $k$ shuffle followed by a cut. Let $N_i$ be the number of
cycles of $w$ of length $i$. Then the $N_i$ are independent and $N_i$
has a negative binomial distribution with parameters $1/i \sum_{d|i}
\mu(d) (k^{i/d}-1)$ and $(u/k)^i$. Consequently, for fixed $k$ as $n
\rightarrow \infty$, the joint distribution of the number of $i$
cycles after a $k$-shuffle followed by a cut converges to independent
negative binomials with parameters $1/i \sum_{d|i} \mu(d) (k^{i/d}-1)$
and $(1/k)^i$. \end{cor}

\begin{proof} Theorem \ref{cyc} and straightforward manipulations
give that

\begin{eqnarray*}
&& 1+ \frac{k-1}{k} \sum_{n \geq 1} \sum_{w \in S_{n}}
\frac{u^n}{nk^{n-1}} {n+k-cd(w)-1 \choose n-1} \prod_i a_i^{n_i(w)}\\
& = & \prod_{i \geq 1}
(1-\frac{a_i u^i}{k^i})^{-1/i \sum_{d|i} \mu(d) (k^{i/d}-1)}.
\end{eqnarray*} Setting all $a_i=1$ gives the equation

\[ 1+\frac{(k-1)u}{k(1-u)} = \prod_{i \geq 1}
(1-\frac{u^i}{k^i})^{-1/i \sum_{d|i} \mu(d) (k^{i/d}-1)}.\] Taking
reciprocals and multiplying by the first equation gives

\begin{eqnarray*}
&& (\frac{1-u}{1-u/k})+ \frac{(k-1)(1-u)}{k-u} \sum_{n \geq 1} \sum_{w \in
S_{n}} \frac{u^n}{nk^{n-1}} {n+k-cd(w)-1 \choose n-1} \prod_i
a_i^{n_i(w)}\\
& = & \prod_{i \geq 1} (\frac{1-\frac{u^i}{k^i}}{1-\frac{a_i u^i}{k^i}})^{1/i
\sum_{d|i} \mu(d) (k^{i/d}-1)}, \end{eqnarray*} proving the first
assertion of the corollary.

	For the second assertion there is a technique simpler than
that in \cite{DMP}. Rearranging the last equation gives that

\begin{eqnarray*}
&& (\frac{1-u}{1-1/k})+ \sum_{n \geq 1} \sum_{w \in
S_{n}} \frac{(1-u) u^n}{nk^{n-1}} {n+k-cd(w)-1 \choose n-1} \prod_i
a_i^{n_i(w)}\\
& = & \frac{1-u/k}{1-1/k} \prod_{i \geq 1} (\frac{1-\frac{u^i}{k^i}}{1-\frac{a_i u^i}{k^i}})^{1/i
\sum_{d|i} \mu(d) (k^{i/d}-1)}. \end{eqnarray*} Letting $g(u)$ be a generating function with a convergent Taylor series, the limit coefficient of $u^n$ in $\frac{g(u)}{1-u}$ is simply $g(1)$. This proves the second assertion. \end{proof}

\section{Open problems} \label{open}

	To finish the paper, we mention two open problems. The most
interesting is to prove the conjecture from \cite{F5} that the cycle
structure of an affine $q$-shuffle is given by the factorization type
of a monic degree $n$ polynomial over $F_q$ {\it with constant term
1}. For the identity conjugacy class, it amounts to the $m=0$ case of
the following observation.

{\bf Corollary (loc. cit.):} For any positive integers $x,y$, the
number of ways (disregarding order and allowing repetition) of writing
$m$ (mod $y$) as the sum of $x$ integers of the set $0,1,\cdots,y-1$
is equal to the number of ways (disregarding order and allowing
repetition) of writing $m$ (mod $x$) as the sum of $y$ integers of the
set $0,1,\cdots,x-1$.

	More generally, let $f_{n,k,d}$ be the coefficient of $z^n$ in
$(\frac{z^k-1}{z-1})^d$ and let $\mu$ be the Moebius function. Let
$n_i(w)$ be the number of $i$-cycles in a permutation $w$. Then
(loc. cit.) the conjecture is equivalent to the truly bizarre
assertion (which we intentionally do not simplify) that for all
$n,k$,

\begin{eqnarray*}
& & \sum_{m=0 \ mod \ n} Coef. \ of \ q^m u^n t^k \
in \ \sum_{n=0}^{\infty} \frac{u^n} {(1-tq)\cdots(1-tq^n)} \sum_{w
\in S_n} t^{cd(w)} q^{maj(w)} \prod x_i^{n_i(w)}\\ & = &
\sum_{m=0 \ mod \ k-1} Coef. \ of \ q^m u^n t^k \ in \
\sum_{k=0}^{\infty} t^k \prod_{i=1}^{\infty} \prod_{m=0}^{\infty}
(\frac{1}{1-q^mx_iu^i})^{1/i \sum_{d|i} \mu(d) f_{m,k,i/d}}.
\end{eqnarray*} This assertion is bizarre because on one side the summation is taken
mod $n$, and on the other side it is taken mod $k-1$!

	A second problem, considered in \cite{F6} is to determine
whether or not the following statement is true. Recall that the major
index of a permutation is the sum of the positions of its descents.

{\bf Statement:} For $n \geq 1$, let $t$ be the largest divisor of $n$
such that $gcd(cd-1,t)=1$. Then for every conjugacy class $C$ of
$S_n$, the set of permutations in $C$ with $cd$ cyclic descents has
its major index equidistributed mod $t$.

\section{Acknowledgements} This research was supported by an NSF
Postdoctoral Fellowship.

\end{document}